\documentclass{amsart}
\usepackage{setspace}
\usepackage[T1]{fontenc}
\usepackage{amssymb,amsmath}
\usepackage{epsfig}
\usepackage{graphicx}
\usepackage{mathtools}
\usepackage{nicefrac}   
\usepackage{lipsum}
\usepackage{amscd}
\usepackage{tikz-cd}
\input xy
\xyoption{all}

\theoremstyle{plain}
\newtheorem{thm}{Theorem}[section]

\theoremstyle{definition} \theoremstyle{definition}

\theoremstyle{remark}

\def\Q{{\mathbb Q}}
\def\Z{{\mathbb Z}}

\def\H{{\mathbb H}}
\def\A{{\mathbb A}}

\def\F{{\overline{F}}}
\def\G{{\mathbb{G}}}

\def\Hom{{\rm Hom}}

\def\Hom{{\rm Hom}}

\def\inv{{\rm inv}}

\title{The Picard Group of a Reductive Group over a Number Field}

\author{Dylon Chow}

\begin{document}

\begin{abstract}
We give a description of the Picard group of a reductive group over a number field as an abelianized Galois cohomology group. It gives another approach of a result due to Labesse.
\end{abstract}

\maketitle

\section{Introduction}

Let $G$ be a connected reductive group defined over a number field $F$. This note is concerned with the Picard group $\text{Pic}(G)$ of $G$. It is known to be finite. Its Pontryagin dual is also a finite group with the same order. It was identified in a paper by Labesse \cite{Labesse1999} as a certain abelian Galois cohomology group. For the notion of abelianized Galois cohomology, we refer the reader to \cite{borovoi1998}. The purpose of this note is to give an alternate, more functorial approach to this bijection.

To state the main result of this paper, we introduce some notation. Let $G^{sc}$ be the simply connected cover of the derived group of $G$. We consider the natural homomorphism $\rho: G^{sc} \rightarrow G$. Let $T$ be a maximal torus in $G$ and let $T^{sc}=\rho^{-1}(T)$. 

Let $\overline{F}$ be an algebraic closure of $F$. For each number field $K$, we let $\A_K$ denote the ring of adeles of $K$. Let $\overline{\A}=\varinjlim \A_K$, the direct limit of the rings $\A_K$, taken over all the finite Galois extensions of $F$ contained in $\overline{F}$. We define the adelic abelian cohomology groups \[\H^i_{ab}(\A_F/F,G):=\H^i(F,T^{sc}(\overline{\A})/T^{sc}(\F) \rightarrow T(\overline{\A})/T(\F)).\] Labesse proved the following result: 

\begin{thm} (Labesse, \cite[Prop. 1.7.3]{Labesse1999})
Let $G$ be a connected reductive group over a global field $F$. There is a canonical bijection  \[H^1_{ab}(\A_F/F,G) \to \text{Hom}(\text{Pic}(G),\Q/\Z).\]
\end{thm}

Labesse proved the result in a quite complicated way. He first used the isomorphism constructed by Kottwitz that described the Picard group of $G$ in terms of the center $Z(G^\vee)$ of the connected component of the Langlands dual group $G^\vee$ \cite{Kot84}. Then he used an extension of the Tate-Nakayama theorem (due to Nyssen) that extends an isomorphism due to Langlands on the Langlands correspondence for tori, also described in terms of $Z(G^\vee)$. 

Our approach avoids using the center of the dual group and formulates our results in terms of the algebraic fundamental group of $G$. This construction is functorial for all morphisms. We use a flasque resolution of reductive groups, constructed by Colliot-Thelene. Similar results on the dual of the Picard group for reductive groups over non-archimedean local fields are known. Our main ingredient is an extension of the Tate-Nakayama theorem to complexes of tori, due to Kottwitz and Shelstad \cite{KS99}.

The proof in this short note may be well-known but I don't know of any reference. It was written in 2021 when I was trying to learn about Galois cohomology by reading some articles of Borovoi and Kottwitz from the 1980's. Before getting into the details, I would like to thank Mikhail Borovoi for answering many of my questions via e-mail and for showing me some papers on flasque resolutions.

\section{Tate-Nakayama Duality for Complexes of Tori}

In this section we review some class field theory. Let $F$ be a number field, $\F$ an algebraic closure, and $\Gamma=\text{Gal}(\F/F)$ the Galois group of $\F$ over $F$. For every finite extension $K$ of $F$ in $\F$, let $\A_K$ be the ring of adeles of $K$, $J_K$ the group of ideles of $K$, and $C_K=J_K/K^\times$ the idele class group.

Let $\overline{\A}_F=\varinjlim \A_K$, the direct limit of the rings $\A_K$ (the limit being taken over the directed set of finite extensions $K$ of $F$ in $\F$. Let $J_{\F}=\overline{\A}_F^\times = \varinjlim \A_K/K^\times=\varinjlim J_K$, and let $C_{\F} = \varinjlim C_K = \overline{\A}_F^\times/\F^\times$.

The group $\overline{\A}_F$ is a smooth $\Gamma$-module, and $\A_K$ can be viewed as the fixed points of $\text{Gal}(K/F)$ in $\overline{\A}$. It is possible to define the set $T(\A_K)$ of points of $T$ with values in $\A_K$, and we let $T(\overline{\A})=\varinjlim T(\A_K)$.

For a linear algebraic group $G$ over a field $F$, we write $G_*$ for the group of cocharacters and $G^*$ for the group of characters.

We have $T(\overline{\A})=T_* \otimes_\Z J_F$, and $T(\F)=T_* \otimes_\Z \F^\times$. A character $\chi \in X^*(T)$ defines compatible maps \[T(\F) \to \F^\times \ \text{and} \ T(\overline{\A}) \to J_F,\] and hence a map \[T(\overline{\A})/T(\F) \to J_F / \F^\times = \varinjlim C_K,\] which induces cup product pairings \[H^r(G_F,X^*(T)) \times H^{2-r}(G_F,T(\overline{\A})/T(\F)) \to H^2(G_F,\varinjlim C_K)=\Q/\Z.\] 
For $r=1$, the cup product pairing induces an isomorphism \[H^1(F,T(\overline{\A})/T(\F)) \to \text{Hom}(H^1(F,X^*(T)),\Q/\Z).\] (see \cite[Lemma D.2.A]{KS99}).

Let $T$ and $U$ be two $F$-tori, and $f \colon T \to U$ a morphism. There is a long exact sequence
\[\dots \to H^i(F,T(\F) \to U(\F)) \to H^i(F,T(\overline{\A}) \to U(\overline{\A})) \to H^i(F,T(\overline{\A})/T(\F)) \to U(\overline{A})/U(\F)) \to \dots\] 
For complexes of length two, the two spectral sequences of hypercohomology reduce to long exact sequences. One of the long exact sequences is \[\dots \to H^r(G,A \to B) \to H^r(G,A) \to H^r(G,B) \to H^{r+1}(G,A \to B) \to \dots \] where the map $i$ is given by \[(a,b) \mapsto a\] for any hypercocycle $(a,b) \in C^r(G,A) \oplus C^{r-1}(G,B)$, the map $f$ is induced by $f \colon A \to B$, and the map $j$ is given by \[b \to (0,b)\] for any cocycle $b$ in $C^r(G,B)$ (see \cite[A.2]{KS99}).

We define \[H^2(\A_F/F,\G_m) = H^2(F,C_{\F})= H^2(\Gamma_F,\overline{\A}_F^\times/\F^\times).\]

Recall that if $K/F$ is a finite Galois extension with group $\Gamma(K/F)$, then for any place $v$ of $F$, we have the decomposition subgroup $\Gamma_v = \Gamma(K_v/F_v)$ associated with $v$. We define \[\inv_{K/F} \colon H^2(\Gamma,I_K) \to \frac{1}{[K:F]}\Z/\Z \subset \Q/\Z\] by the formula \[\inv_{K/F}(c) = \sum_{v \in \Omega_F}\inv_{K_v/F_v}(c_v),\] where $c_v$ is the component at $v$ of $c \in H^2(\Gamma(K/F),I_K)=\bigoplus_{v \in \Omega_F}H^2(\Gamma_v,K_v^\times).$

Passing to the limit over the finite Galois extensions $K$ of $F$, we obtain a map \[\inv_F \colon H^2(F,I) \to \Q/\Z,\] which induces an isomorphism \[\inv_F \colon H^2(\A_F/F,\G_m)=H^2(F,C) \to \Q/\Z.\]

We define groups \[\H^i(\A/F,T \to U):=\H^i(F,T(\overline{\A})/T(\F) \to U(\overline{\A})/U(\F)).\]
There is a Tate-Nakayama pairing (\cite[C.2.1]{KS99}) \[\H^r(\A_F/F,T \to U) \otimes \H^{2-r}(F,U^* \to T^*) \to \Q/\Z,\] the $\Q/\Z$ coming from the canonical isomorphism \[\inv_F \colon H^2(\A_F/F,\G_m)=H^2(F,C) \to \Q/\Z.\]
Apply the exact sequence for $T \to U$ and $U^* \to T^*$. The 5-lemma implies that the Tate-Nakayama pairing induces an isomorphism \[H^1(\A/F,T \to U) \to \Hom(H^1(F, U^*\to T^*),\Q/\Z)\] (\cite[Lemma C.2.A]{KS99}). 

\section{Picard Group as an Abelian Cohomology Group}

In this section we explain how the Picard group of $G$ may be viewed as a hypercohomology group. Such an isomorphism was constructed by Borovoi and van Hamel in a general setting (\cite{BvH2009}), but there is an easier construction when the variety is an algebraic group. This may be well-known but I cannot find this result in the literature. I am grateful to Mikhail Borovoi for indicating a proof. 

By \cite[Proposition-Definition 3.1]{Colliot2008}, there exists a resolution \[1 \to S \to H \to G \to 1,\] where $S$ is an $F$-torus and $H$ is a quasi-trivial reductive group, that is, $[H,H]$ is simply connected and $P:=H/[H,H]$ is a quasi-trivial $F$-torus. 

By \cite[Proposition 3.3]{Colliot2008}, we have an exact sequence \[(P^*)^\Gamma \to (S^*)^\Gamma \to \text{Pic} \ G \to 0.\] 
The hypercohomology group $\H^1(F,P^* \to S^*)$ fits into the exact sequence \[(P^*)^\Gamma \to (S^*)^\Gamma \to \H^1(F,P^* \to S^*) \to H^1(F,P^*) = 0.\] Thus we obtain an isomorphism \[\text{Pic} \ G \cong \H^1(F,P^* \to S^*).\] There is a canonical isomorphism \[\H^1(F,T^* \to T_{sc}^*) \cong \H^1(F,P^* \to S^*),\] which gives an isomorphism \[\text{Pic}(G) \cong \H^1(F,T^* \to T_{sc}^*).\]
Applying this we get \[\Hom(\text{Pic}(G),\Q/\Z) = \Hom(\H^1(F,T^* \to T_{sc}^*),\Q/\Z) = H^1(\A/F,T^{sc} \to T) =:H^1_{ab}(\A_F/F,G).\] This completes the proof.

\bibliographystyle{alpha}
\bibliography{biblio.bib}

\end{document}